\documentclass[12pt]{article}
\usepackage{amsfonts,amssymb,amscd}
\usepackage{oldgerm}

\usepackage[Conny]{fncychap}
\usepackage[latin1]{inputenc}
\usepackage{amssymb,amsmath,amsthm,amsfonts}
\usepackage[latin1]{inputenc}
\usepackage{mathrsfs}
\usepackage[all]{xy}

\newcommand{\ber}{\mathcal{B}}

\newtheorem{teo}{Theorem}[section]

\newtheorem{lema}[teo]{Lemma}

\newtheorem{defnc}[teo]{Definition}

\newtheorem{obs}[teo]{Remark}

\newcommand{\R}{{\mathbb R}}

\newcommand{\nat}{\mathbb{N}}

\newcommand{\rar}{\rightarrow}
\begin{document}

\title{The Analyticity of a Generalized Ruelle's Operator}


\author{R.R.Silva, E. A. Silva, R. R. Souza}
\date{\today}
\maketitle

\begin{abstract}
In this work we propose a generalization of the concept of Ruelle's operator for one dimensional lattices used in thermodynamic formalism and ergodic optimization, which we call generalized Ruelle's operator, that generalizes both the Ruelle  operator proposed in \cite{le} and the Perron Frobenius operator defined in \cite{bo}. We suppose the alphabet is given by a
compact metric space, and consider a general a-priori measure to define the operator. We also consider the case where the set of symbols that can follow a given symbol of the alphabet
depends on such symbol, which is an extension of the original concept of transition matrices from the theory of subshifts of finite type.
 We prove the analyticity of the Ruelle´s operator
  and present some examples.
\end{abstract}

\bigskip


\bigskip


\section{Introduction}

\hspace{0.4cm}

The Ruelle´s Operator (or transfer operator) has been used since \cite{bo,Ruelle} to describe equilibrium states in statistical mechanics, in what is called the classical thermodynamical formalism on the Bernoulli space $M^{\nat}$, where $M$ is called the {\it alphabet}. In such setting one usually considers the set $M$ as a finite set (\cite{PP}). Using Markov partitions, one can prove that the shift map on the Bernoulli space can be used to describe hyperbolic maps, in what is known as the symbolic dynamics setting.
The thermodynamical formalism can also be used to study invariant measures which have regular Jacobians (see \cite{Sar2}).

We fix a potential $\psi:M^{\nat} \rar \R$ that describes an interaction scheme in statistical mechanics. Usually such potential is Holder-continuous. We then define the Ruelle´s operator on a certain subset of continuous functions (usually Holder) from $M^{\nat}$ to $\R$ as $$L_{\psi}(\varphi)(x)=\sum_{a\in M} e^{\psi(ax)}\varphi(ax),$$ where $x=(x_1,x_2,x_3,\hdots) \in M^{\nat}$ and $ax=(a,x_1,x_2,x_3,\hdots)$.
The eigenfunctions and eigenmeasures associated to the maximal positive eigenvalue of the Ruelle´s operator can be used to define an invariant and ergodic measure (called equilibrium measure) that satisfies a variational principle, maximizing the sum of the integral of the potential $\psi$ with the entropy of the measure (or minimizing {\it free energy}), where the maximization is performed among all  the invariant measures for the shift map
(see also \cite{Keller}). Under some conditions, when we multiply the potential $\psi$ by a constant $\beta$ (the inverse of the temperature), and send $\beta \rightarrow 0$, the equilibrium measure  converge (in weak sense) to the invariant measure that maximizes the integral of the potential $\psi$. This last problem is usually known as an Ergodic Optimization problem, or the analysis of the zero temperature case.
(see \cite{Bousch, Bremont, CLT, Jenkinson, Mane}).

An interesting problem is to analyze what happens when we have a different alphabet. For example \cite{Sar} proposed models where the alphabet $M$ is an infinite countable space, such as $\nat$. This setting can be used to describe non-uniformly hiperbolic maps, for example, the Manneville-Pomeau maps.
Questions of Ergodic optimization in this setting where also considered (see \cite{BGar, JMUrb,Sar2,Daon}).

In the present paper, as proposed in \cite{LMST}, we consider the case where $M$ is a compact metric space. For example, when  $M=S^1$, we have the so called XY-model.


In the paper \cite{LMST} the thermodynamical formalism for potentials that depends on two coordinates (see \cite{spit}) on a compact metric space (which is the case of the XY model) was considered. In the XY spin model (see \cite{FH} and references therein), we can have for example a potential $\psi(x, y) = cos(x.y)$, where $x, y \in (0, 2\pi]$. When there is a magnetic term one could consider, for instance, $\psi(x, y) = cos(x . y) + m cos(x)$, where $m$ is constant \cite{A}.
In another paper \cite{le}, several thermodynamic formalism questions where analyzed for the case where the potential can depend on all the coordinates of the Bernoulli space.
In both works, Ruelle´s operator can be used to prove the existence of equilibrium states, and also consider ergodic optimization problem, where we search for the invariant measure that maximizes the integral of the potential. One can prove that, under uniqueness assumptions for the maximizing measure, such measure is the weak limit of the equilibrium states when the temperature decreases to zero. \cite{LMST} and \cite{LM} obtain large deviation principles for this convergence, and also questions involving {\it selection} of measures (see also \cite{vanenter, CH}).

As pointed out in \cite{le}, a major difference between the settings of the classical Bowen-Ruelle-Sinai 
Thermodynamic Formalism and the XY model is that here, in order to
define the Ruelle´s operator, we need an a-priori measure defined on $M$.
In \cite{LMST,le} the a-priori measure was considered to be the Lebesgue measure. 
In the present paper, we analyse the Ruelle's Operator 
 where the a-priori measure is any measure on the Borel sets of $M$,
 and  also propose a formulation of the Ruelle Operator, which we call the Generalized Ruelle's Operator, that generalizes both the Ruelle  Operator proposed in \cite{le} as well as the Perron Frobenius Operator defined in \cite{bo}, by  considering the case where there are restrictions on  transitions (i.e. a generalization of the concept of subshift - see section \ref{subshift}).
We also prove the analyticity of this Generalized Ruelle's Operator, and present some Taylor
expansions of the operator  that may be useful in future works (see for instance equation \eqref{pse}).


Now, we describe in brief what was developed during this work: Initially, let us fix $M$ a compact metric space and consider the set ${\cal B}$ of all sequences of points in $M$ (paths).
We shall denote by $C^{\gamma }({\cal B}, \mathbb R )$ the usual  H\"older Space. For each potential $\psi $ belonging to $C^{\gamma }({\cal B}, \mathbb R )$ we define the  Ruelle Operator ${\mathscr L}_{\psi}:C^{\gamma }({\cal B}, \mathbb R )\hookleftarrow$ associated to $\psi$
 $$
{\mathscr{L}}_{\psi}(\varphi )(x)=\int \limits_M e^{\psi (ax)} \varphi(ax) d\mu (a)\;,\;\;\varphi \in C^{\gamma},\; x \in \cal B\,,
$$
where $ax$ denotes the sequence $ax=(a,x_0,x_1,\ldots)$.
We prove the analyticity of this operator with respect to $\psi$.

Then we replace ${\cal B}$ by ${\cal B}(A,I)$, where $I$ is a closed subset of the real line and $A:M\times M \rightarrow \mathbb R$ is a sectionally trivial function on $I$ (see Definition 2).
 Note that $A$ and $I$ play the role of the transition matrices in the theory of subshifts of finite type.
We present some examples that show that the classical setting of subshifts of finite type is indeed a particular case of our setting.
Then we will present an example which was introduced in \cite{LMMS} and can be used to understand the countable alphabet case.
For every $\psi $ belonging to $C^{\gamma }({\cal B}(A,I), \mathbb R )$ we define the Generalized Ruelle's Operator
${\mathscr L}_{\psi}:C^{\gamma }({{\cal B}(A,I)}, \mathbb R )\hookleftarrow$ associated to $\psi $ (see Definition~(\ref{sruelle})).
We prove the analyticity of this operator with respect $\psi$, which implies the analyticity of the Ruelle operator in the examples given before.



\bigskip


\section{Preliminaries}

\hspace{0.5cm}
In this section we establish some preliminary results and definitions which
will be used throughout this work.

Consider $M=(M,d)$ a compact metric space, equipped with a certain finite measure $\mu$ defined on the borelian sets of $M$.

We shall denote by $\cal B$ the set of all sequences $x=(x_0,x_1,\ldots )$, $ x_i \in M$,  $i\in {\mathbb N}$, equipped with the distance $d_c$, $c>1$, given by 

$$d_c (x,y)=\sum \limits_{k\geq 0} \frac{d(x_k ,y_k )}{c^k}\; , \;\;x,y \in \cal B$$

The topology generated by $d_c$ is the product topology. It follows from Tychnoff's Theorem that $\cal B$ is a compact metric space.

We shall also be concerned with the space $C^{\gamma}$ equipped with the norm $\mid . \mid _\gamma$, $0\leq \gamma \leq 1$, defined as follows:

- in the case $\gamma =0$, we denote by $C^0= C^0(\cal B, \mathbb R )$  the set of all continuous functions from $\cal B$ into $\mathbb R$ and define $ \mid \varphi \mid  _0 = \sup \limits_{x\in \cal B} |\varphi (x) | \; , \;\;  \forall \varphi \in C^{0}$.

- in the case  $0<\gamma \leq 1$, we denote by  $C^{\gamma}$ the set of all $\gamma $- H\"older continuous functions from $\cal B$ into $\mathbb R$, i.e.,
$$C^{\gamma}=C^{\gamma} ({\cal B},{\mathbb R})=\left\{ \varphi:{\cal B} \rightarrow {\mathbb R}\; ; \;\; \sup \limits_{x \neq y} \frac{\mid \varphi (y) -\varphi (x) \mid}{d_c ^{\gamma}(x,y)}<\infty \right\}$$
and we define the norm $\mid \varphi \mid _{\gamma} = \mid \varphi \mid _0 + Hol_{\varphi} \; , \;\; \forall \varphi \in C^{\gamma}$, where $Hol_{\varphi}$ is given by

$$Hol_{\varphi}= \sup \limits_{x \neq y} \frac{\mid \varphi (y) -\varphi (x) \mid}{d_c ^{\gamma}(x,y)}\; , \;\;\forall \varphi \in C^{\gamma}.$$

It is well known that $C^{\gamma}$ is a Banach space, $\forall \,\,0\leq \gamma \leq1$.

\begin{obs}
\label{holdernorm}

We list below some elementary properties of the norm $\mid .\mid _{\gamma}$:
\begin{enumerate}
\item $\mid \varphi \psi \mid _{\gamma}\leq  2 \mid \varphi \mid _{\gamma} \mid \psi \mid _{\gamma}$, for all $\varphi$ and $\psi$ in $C^{\gamma}$.
\item $\mid \beta _1  \ldots \beta _k \mid _{\gamma} \leq 2^{k-1} \mid \beta _1 \mid _{\gamma} \ldots  \mid \beta _k \mid _{\gamma}$, for all $ \beta _1, \ldots , \beta _k$ in $C^{\gamma}$.
\end{enumerate}

\end{obs}

\noindent {\bf{ Proof}}.
In fact, property (2) is an immediate consequence of property (1).
To prove item (1), initially we note that for every $x$ and $y$ in $\cal B$, we have,
\begin{align*}
|(\psi \varphi)(x) - \psi \varphi)(y)| \leq  &
|\psi (x) \varphi (x) - \psi (x) \varphi (y)|+ |\psi (x) \varphi (y) - \psi (y) \varphi (y)|  \\
\leq & | \psi (x) | | \varphi (x) - \varphi (y) | + | \varphi (y) ||\psi (x) - \psi (y)| \\
\leq  & \mid \psi \mid _0 Hol _{\varphi} d_c ^{\gamma}(x,y) + \mid \varphi \mid _0 Hol_{\psi} d_c^{\gamma}(x,y)\,.
\end{align*}
Thus, if $x\neq y$ we have

$$\frac{|(\psi \varphi)(x) - (\psi \varphi )(y)|}{d_c ^{\gamma}(x,y)}\leq \mid \psi \mid _0 Hol _{\varphi} + \mid \varphi \mid _0 Hol_{\psi}
\leq 2 \mid \psi \mid _{\gamma}\mid \varphi \mid _{\gamma}\,. $$
The inequality above implies that $Hol_{\psi \varphi}\leq 2 \mid \psi \mid _{\gamma}\mid \varphi \mid _{\gamma}$.
Since
$\mid \psi \varphi \mid _0  \leq  \mid \psi \mid _{0} \mid \varphi \mid _{0} \leq \mid \psi \mid _{\gamma} \mid \varphi \mid _{\gamma}$, we may conclude $(1)$.
\qed


\section{The Ruelle's Operator}

\subsection{Definition and basic properties}

\begin{defnc}
\label{ruelle}
For any $\psi \in C^{\gamma}, \; 0\leq \gamma \leq 1$, the Ruelle's Operator ${\mathscr{L}}_{\psi}:C^{\gamma}\hookleftarrow$ associated to $\psi$ is given by
$$
{\mathscr{L}}_{\psi}(\varphi )(x)=\int \limits_M e^{\psi (ax)} \varphi(ax) d\mu (a)\;,\;\;\varphi \in C^{\gamma},\; x \in \cal B\,,
$$
where $ax$ denotes the sequence $ax=(a,x_0,x_1,\ldots)$, for any $a\in M$ and $x\in \cal B$.

\end{defnc}

\begin{obs}
\label{bemdefinida}
For any $\psi \in C^{\gamma}, \; 0\leq \gamma \leq 1$, the Ruelle's Operator ${\mathscr{L}}_{\psi}:C^{\gamma}\hookleftarrow$ associated to $\psi$ is well defined and it is a bounded linear operator. More precisely, we have

\begin{enumerate}

\item  If $\psi \in C^0$ then ${\mathscr{L}}_{\psi}(\varphi ) \in C^0$, for all $\varphi \in C^0$.

\item If $\psi \in C^{\gamma},\; 0<\gamma \leq 1$, then ${\mathscr{L}}_{\psi}(\varphi ) \in C^{\gamma}$, for all $\varphi \in C^{\gamma}$.

\item  If $\psi \in C^0$ then ${\mathscr{L}}_{\psi}:C^0\hookleftarrow$ is a bounded linear operator.

\item  If $\psi \in C^{\gamma},\; 0<\gamma \leq 1$, then  ${\mathscr{L}}_{\psi}:C^{\gamma}\hookleftarrow$ is a bounded linear operator.

\end{enumerate}
\end{obs}

\noindent {\bf{ Proof}}. 

{\bf {1.}} For every $x$ and $y$ in $\cal B$, we have,
\begin{align*}
|{\mathscr{L}}_{\psi}(\varphi )(x)-{\mathscr{L}}_{\psi}(\varphi )(y)|  = & |\int \limits_M e^{\psi (ax)} \varphi (ax) d\mu (a)- \int \limits_M e^{\psi (ay)} \varphi (ay) d\mu (a)|  \\
 \leq & \int \limits_M |e^{\psi (ax)} \varphi (ax) - e^{\psi (ay)} \varphi (ay)| d\mu (a)\,.
\end{align*}

 At this point we need to estimate the integrand above.
\begin{align*}
|e^{\psi (ax)} \varphi (ax) - e^{\psi (ay)} \varphi (ay)|  = & |e^{\psi (ax)} \varphi (ax) - e^{\psi (ax)} \varphi (ay)+e^{\psi (ax)} \varphi (ay) - e^{\psi (ay)} \varphi (ay)|  \\
 \leq & |e^{\psi (ax)} (\varphi (ax) - \varphi (ay))|+|\varphi (ay)(e^{\psi (ax)} - e^{\psi (ay)})|\,.
\end{align*}

By the Mean Value Theorem, for every $z$ and $w$ in ${\mathbb R}$, we have $e^z -e^w =e^{\hat {z}} (z-w)$, for some $\hat {z}$ between $z$ and $w$. Thus, we have,
$|e^{\psi (ax)}  - e^{\psi (ay)} |\leq e^{|\psi|_0} |\psi (ax) - \psi (ay)|$,
and 
therefore,
\begin{equation}\label{integrando}
    |e^{\psi (ax)} \varphi (ax) - e^{\psi (ay)} \varphi (ay)|\leq  e^{|\psi|_0} \Big(|\varphi (ax) - \varphi (ay)|+ \mid \varphi \mid _0 |\psi (ax) - \psi (ay)| \Big)\,.
\end{equation}

Now, fix $\varepsilon >0$. Since $\cal B$ is compact $\varphi $ and $\psi$ are uniformly continuous functions. Thus, there exists $\delta _1>0$ such that
$$
\max \left\{|\varphi (t) - \varphi (s)|,\; |\psi (t) -\psi (s)|\right\}<\frac{\varepsilon}{\mu (M)e^{|\psi|_0}  (1+\mid \varphi \mid _0)}, \mbox{ if }d_c(t,s)<\delta _1,\; t,s \in \cal B\,.
$$

Using $d_c(ax,ay)=\frac{d(x,y)}{c}$ and \eqref{integrando}, if $d_c(x,y)< c\delta_1$ we have 

 $$|{\mathscr{L}}_{\psi}(\varphi )(x)-{\mathscr{L}}_{\psi}(\varphi )(y)|\leq \int \limits_M e^{|\psi|_0} \left\{\frac{\varepsilon}{\mu (M)e^{|\psi|_0} (1+\mid \varphi \mid _0)} + \mid \varphi \mid _0\frac{\varepsilon}{\mu (M)e^{|\psi|_0} (1+\mid \varphi \mid _0)}\right\}d\mu (a) = \varepsilon,$$
which implies that ${\mathscr{L}}_{\psi}(\varphi ):\cal B \rightarrow \mathbb R$ is a uniformly continuous function, and in particular  ${\mathscr{L}}_{\psi}(\varphi ) \in C^0$.

{\bf {2.}} From \eqref{integrando} above we have
$$|e^{\psi (ax)} \varphi (ax) - e^{\psi (ay)} \varphi (ay)|
\leq
e^{|\psi|_0} \left\{ Hol_{\varphi } d_c ^{\gamma} (ax,ay) + \mid \varphi \mid _0 Hol _{\psi}d_c^{\gamma}(ax,ay)           \right\}\,.$$

We observe that $d_c(ax,ay)=\frac{d(x,y)}{c}$, and thus we have

\begin{align}
\nonumber
|{\mathscr{L}}_{\psi}(\varphi )(x)-{\mathscr{L}}_{\psi}(\varphi )(y)|\leq
&
\frac{1}{c^{\gamma}} \int \limits_M e^{|\psi|_0} \left\{ Hol_{\varphi } d_c ^{\gamma} (x,y) + \mid \varphi \mid _0 Hol _{\psi}d_c^{\gamma}(x,y) \right\} d\mu(a)
\\
 = & \frac{\mu (M)e^{|\psi|_0}}{c^{\gamma}}\left\{ Hol_{\varphi} + \mid \varphi \mid _0 Hol _{\psi}\right\} d_c ^{\gamma}(x,y)\,. \label{holoper}
\end{align}

Hence
$$
\sup \limits_{x\neq y}\frac{|{\mathscr{L}}_{\psi}(\varphi )(x)-{\mathscr{L}}_{\psi}(\varphi )(y)|}{d_c ^{\gamma}(x,y)}
\leq
\frac{\mu (M)e^{|\psi|_0}}{c^{\gamma}}\left\{ Hol_{\varphi} + \mid \varphi \mid _0 Hol _{\psi}\right\},
$$
which means that ${\mathscr{L}}_{\psi}(\varphi ) \in C^{\gamma}$.

{\bf {3.}} Clearly, ${\mathscr{L}}_{\psi}:C^{\gamma}\hookleftarrow$ is linear operator.
Now, given  $\varphi \in C^0$ we have,
$$\mid {\mathscr{L}}_{\psi}(\varphi ) \mid _0 =
\sup\limits_{x\in \cal B}|\int \limits_M e^{\psi (ax)} \varphi (ax)d\mu (a)| \leq
\mu (M) e^{\mid \psi \mid _0} \mid \varphi \mid _0\,.$$
Therefore,
${\mathscr{L}}_{\psi}:C^0\hookleftarrow$ is a bounded linear operator.

{\bf {4.}}  Here we also have that ${\mathscr{L}}_{\psi}:C^0\hookleftarrow$ is a linear operator.
It follows from \eqref{holoper} above that for any $\varphi \in C^{\gamma}$
we have,
\begin{align*}
Hol_{{\mathscr{L}}_{\psi}(\varphi )}\leq
& \frac{\mu (M)e^{|\psi|_0}}{c^{\gamma}} \left\{ Hol_{\varphi} + \mid \varphi \mid _0 Hol _{\psi}\right\}  \\
 \leq & \frac{\mu (M)e^{|\psi|_0}}{c^{\gamma}} \left\{ \mid {\varphi} \mid _{\gamma} + \mid \varphi \mid _{\gamma} \mid {\psi} \mid _{\gamma}\right\} \\
 = & \frac{\mu (M)e^{|\psi|_0}}{c^{\gamma}}\Big( 1 + \mid {\psi} \mid _{\gamma}\Big)\mid {\varphi} \mid _{\gamma}\,.
\end{align*}

We also have
 $$\mid {\mathscr{L}}_{\psi}(\varphi ) \mid _0 \leq  \mu (M) e^{\mid \psi \mid _0} \mid \varphi \mid _0 \leq \mu (M) e^{\mid \psi \mid _0} \mid \varphi \mid _{\gamma},$$
where the first inequality was proved in the proof of item ({\it {3}}).

Hence
$$\mid {\mathscr{L}}_{\psi}(\varphi ) \mid _{\gamma} =
\mid {\mathscr{L}}_{\psi}(\varphi ) \mid _0 + Hol_{{\mathscr{L}}_{\psi}(\varphi )} \leq
\left\{ \frac{\mu (M)e^{|\psi|_0}}{c^{\gamma}}\Big(1+\mid \psi \mid _{\gamma}\Big)+ \mu (M) e^{\mid \psi \mid _0} \right\}\mid {\varphi} \mid _{\gamma}$$
which implies that  ${\mathscr{L}}_{\psi}:C^{\gamma}\hookleftarrow$ is a bounded linear operator.

\qed

\vspace{1cm}


\subsection{Analyticity}

The main result of this section is theorem \ref{principal},
which states that the function $\Theta :C^{\gamma}\rightarrow V$ that sends $\psi \in C^{\gamma}$ to the
Ruelle's Operator ${\mathscr{L}}_{\psi} \in V$ associated to the potential $\psi$ (
$V$ is the set of all continuous linear transformations $l:C^{\gamma}\hookleftarrow$
), is an analytic function.
 Before proving such result, we have to discuss the concept of analyticity we are dealing with in this paper.
 After proving theorem \ref{principal}, we show a Taylor expansion of the Ruelle´s operator that may be useful in applications.

\begin{defnc}
\label{der-ordemsup}
Let $(X, \|\cdot \|_X)$ and $(Y, \|\cdot \|_Y)$ be Banach spaces and $U$ an open subset of $X$. 
If $k\in \mathbb{N}$, a function $F:U\rightarrow Y$ is said to be $k$-differentiable in $x$ if, for each $j \in \{1,...,k\}$, there exist a $j$-linear bounded transformation $D^j F(x): X^j \rightarrow Y$ (where $X^j$ is the product space given by $j$ copies of $X$) such that
$$
D^{j-1}F(x+v_j)(v_1,...,v_{j-1}) - D^{j-1}F(x)(v_1,...,v_{j-1}) = D^j F(x) (v_1,...,v_j) + o_j(v_j)
$$
where $o_j:X \rightarrow Y$ is such that $\lim_{v\rar 0} \frac{\|o_j(v)\|_Y}{\| v \|_X} = 0$
\end{defnc}

We say that $F$ has derivatives of all orders in $U$, if for any $k\in \mathbb{N}$, and any $x \in U$, F is $k$-differentiable in $x$.

\begin{defnc}
\label{analiticidade}

Let $X$ and $Y$ be Banach spaces and $U$ an open subset of $X$. A function $F:U\rightarrow Y$ is called analytic on $U$, when $F$ has derivatives of all orders in $U$, and for each $x\in U$ there exists an open neighborhood $U_x$ of $x$ in $U$ such that for all $v\in U_x$, we have that
$$
F(x+v)-F(x)=\sum \limits_{j=1}^{\infty}\frac{1}{n!}D^jF(x) v^j\,,
$$
where $D^j F(x) v^j=D^jF(x)(v,\ldots ,v)$ and $D^j F(x)$ is the $j-$th derivative of $F$  in $x$.

\end{defnc}

Note that when a function $F:U\rightarrow Y$ is analytic on $U$, then for each $n\in \mathbb{N}$, we have the Taylor expansion of order $n$:
\begin{equation}\label{Taylor}
    F(x+v)=F(x)+D^1F(x)v+\frac{D^2F(x)v^2}{2}+\frac{D^3F(x)v^3}{6}+\hdots + \frac{D^nF(x)v^n}{n!}+o_{n+1}(v)\,,
\end{equation}
where
$o_{n+1}(v)=\sum \limits_{j=n+1}^{\infty}\frac{1}{n!}D^jF(x) v^j$ satisfies $\lim_{v \rar 0} \frac{\|o_{n+1}(v)\|_Y}{\|v\|_X^n}=0$.


\bigskip

In what follows,
the role of $X$ will be played by $C^{\gamma}$, while the role of $Y$ will be played by the set $V$ of all continuous linear transformations $l:C^{\gamma}\hookleftarrow$, endowed with the usual norm
$$
\mid l \mid _{V} =\sup \limits_{\mid \varphi \mid _{\gamma}\leq 1} \mid l( \varphi )\mid _{\gamma}\,,
$$
which makes $V$ a Banach space. 

Now we state and prove our main result.

\begin{teo}
\label{principal}
The function $\Theta :C^{\gamma}\rightarrow V$ that sends $\psi \in C^{\gamma}$ to the
Ruelle's Operator ${\mathscr{L}}_{\psi}$ associated to the potential $\psi$, is an analytic function.
\end{teo}

Before proving this result, we point out that \cite{Sil} presented the proof of analiticyty of the Ruelle operator in the simpler case where we have a finite alphabet and the a-priori measure is the counting measure. See also \cite{BCV} for related results in the finite alphabet case.

\noindent {\bf{ Proof.}}  Given $\psi,\;\beta$ and $\varphi$ in $C^{\gamma}$, for any $x\in \cal B$ we have,

\begin{align*}
\Theta (\psi + \beta)(\varphi)(x) - \Theta (\psi)(\varphi)(x)
= & \;\; {\mathscr{L}}_{\psi +\beta}(\varphi)(x) -{\mathscr{L}}_{\psi}(\varphi)(x)  \\
 = & \int \limits_{M} e^{\psi (ax) +\beta (ax)} \varphi (ax) d\mu (a)-  \int \limits_{M} e^{\psi (ax) } \varphi (ax) d\mu (a) \\
 = & \int \limits_M e^{\psi (ax) } \varphi (ax) \big( e^{\beta (ax)} -1 \big)d\mu (a) \\
 = & \int \limits_M \biggl(e^{\psi (ax) } \varphi (ax) \sum \limits_{n=1} \limits^{\infty} \frac{[\beta (ax)]^n}{n!} \biggr)d\mu (a)\,.
\end{align*}

By Fubini's Theorem 
(we will soon prove that we can use such theorem)
we have,

$$\Theta (\psi + \beta)(\varphi)(x) - \Theta (\psi)(\varphi)(x) =\sum \limits_{n=1} \limits^{\infty}\frac{1}{n!}\int \limits_M e^{\psi (ax)}\varphi (ax)[\beta (ax)]^n d\mu (a)\,.$$

Thus
$$\Theta (\psi + \beta)(\varphi)(x) - \Theta (\psi)(\varphi)(x) = \sum \limits_{n=1} \limits^{\infty}\frac{1}{n!}\Theta (\psi)(\varphi [\beta]^n)(x) $$

If we omit the dependence on $x$ and  $\varphi $ we obtain

\begin{eqnarray}
\label{representacaoemseriesdepotencias}
\Theta (\psi + \beta) - \Theta (\psi) = \sum \limits_{n=1} \limits^{\infty}\frac{1}{n!}\Theta (\psi)((\;.\;) [\beta]^n)\,.
\end{eqnarray}

In order to apply Fubini´s Theorem we have to prove that the sum above converges. To do that, we note that for any $\beta _1 ,\ldots ,\beta _k$ in $C^{\gamma}$, the map
$$
\varphi \mapsto \Theta (\psi )(\varphi \beta _1 \ldots \beta_k )
$$
from $C^{\gamma}$ into $C^{\gamma}$ is a continuous linear map, i.e.,
$\Theta (\psi )((.) \beta _1 \ldots \beta_k ) \in V$. Indeed, we can prove that

\begin{eqnarray}
\label{dominando}
\mid \Theta (\psi )((\;.\;) \beta _1 \ldots \beta_k ) \mid _V \leq \mid {\mathscr L}_{\psi} \mid _V 2^k \mid \beta _1 \mid _{\gamma}\ldots \mid \beta _k \mid _{\gamma}
\end{eqnarray}
because
\begin{align*}
\mid \Theta (\psi )((\;.\;) \beta _1 \ldots \beta_k ) \mid _V = &
\sup _{\mid \varphi \mid _{\gamma}\leq 1}\mid \Theta (\psi )((\varphi) \beta _1 \ldots \beta_k ) \mid _{\gamma}  \\
 \leq & \sup _{\mid \varphi \mid _{\gamma}\leq 1}\mid \mathscr L _{\psi }\mid _V \mid \varphi \beta _1 \ldots \beta _k \mid _{\gamma }\\
 \leq &
 \sup _{\mid \varphi \mid _{\gamma}\leq 1} \mid \mathscr L _{\psi }\mid _V  2^k\mid \varphi  \mid _{\gamma} \mid \beta _1 \mid _{\gamma} \ldots  \mid \beta _k \mid _{\gamma}
 \,,
\end{align*}
where in the last inequality we used ({\it{2} }) of Remark~(\ref{holdernorm}).


Therefore, it follows from inequality~(\ref{dominando}) that
$$
\sum \limits_{n=1} \limits^{\infty}\frac{1}{n!} \mid \Theta (\psi)((\;.\;) [\beta]^n) \mid _V\leq \sum \limits_{n=1} \limits^{\infty}\frac{1}{n!}\mid \mathscr L _{\psi }\mid _V (2\mid \beta \mid _{\gamma})^n = \mid \mathscr L _{\psi }\mid _V \left\{ e^{2 \mid \beta \mid _{\gamma}} -1 \right\}
$$
which implies that the series $\sum \limits_{n=1} \limits^{\infty}\frac{1}{n!}\Theta (\psi)((\;.\;) [\beta]^n)
$ converges in $V$.

Now, the analyticity of $\Theta $ follows at once from the following statement:

\noindent {\bf Statement:} $D^k \Theta (\psi )(\beta _1, \ldots , \beta _k)= \Theta (\psi ) ((\;.\;)\beta _1 ,\ldots , \beta _k)\;,\;\; \beta _1 ,\ldots, \beta_k$ in $C^{\gamma}$, $k=1,2,\ldots $.

\noindent {\bf{ Proof of the Statement.}} The verification will be carried out by induction on $k$.

In what follows  $L^k=L^k(C^{\gamma},V)$ denotes the set of all continuous $k-$linear  functions $l:C^{\gamma} \times \ldots \times C^{\gamma}\rightarrow V$, from $C^{\gamma} \times \ldots \times C^{\gamma}$ ($k-$ copies of   $C^{\gamma } $) into $V$. The norm $\mid .\mid _{L^k}$ of $L^k$, is given by

$$ \mid l \mid _{L^k}=\sup \limits_{\stackrel{\beta _i}{i=1,\ldots ,k}}\mid l\mid _V\;,\;\; l\in L^k\,.$$

Let us prove that the statement is true for $k=1$: in fact, by equation~(\ref{representacaoemseriesdepotencias}) we have,

$$\Theta (\psi + \beta _1) - \Theta (\psi ) = \Theta (\psi )( (\;.\;)\beta _1) + E_1(\beta _1)$$

Where, $E_1(\beta_1 )=\sum \limits_{n=2}\limits^{\infty} \frac{1}{n!}\Theta (\psi )((\;.\;)[\beta_1]^n)$.

It follows from inequality~(\ref{dominando}) that

$$\mid \Theta (\psi )((\;.\;)\beta _1) \mid _V \leq \mid \mathscr L _{\psi} \mid _V 2\mid \beta_1 \mid _{\gamma}$$

Thus, the map $\beta _1 \mapsto \Theta (\psi )((\;.\;)\beta _1)$ is in $L^1$. Again, in view of  inequality~(\ref{dominando}) we have

\begin{align*}
\mid E_1 (\beta _1) \mid _V
= & \mid \sum \limits_{n=2}\limits^{\infty} \frac{1}{n!}\Theta (\psi )((\;.\;)[\beta_1]^n) \mid _V   \\
 \leq & \sum \limits_{n=2}\limits^{\infty} \frac{1}{n!}\mid \mathscr L _{\psi} \mid _V (2\mid \beta \mid _{\gamma})^n  \\
 = & \mid \mathscr L _{\psi} \mid _V \Bigl(e^{2\mid \beta _1 \mid _{\gamma}}-2\mid \beta _1 \mid _{\gamma}-1 \Bigr)
\end{align*}


Hence, $\frac{\mid E_1 (\beta _1)\mid _V}{\mid \beta _1 \mid _{\gamma}}\rightarrow 0$ if $\mid \beta _1 \mid _{\gamma} \rightarrow 0$. Therefore, $D^1 \Theta (\psi)(\beta _1)=\Theta (\psi )((\;.\;)\beta _1)$, and the statement is true for $k=1$.

Now, let us suppose the statement is true for $k-1,\;k\geq 2$, i.e.,

\begin{eqnarray}
\label{k-1-derivada}
D^{k-1} \Theta (\psi )(\beta _1, \ldots , \beta _{k-1})= \Theta (\psi ) ((\;.\;)\beta _1 \ldots \beta _{k-1})\;,\;\; \beta _1 ,\ldots, \beta_{k-1}\in C^{\gamma}\,.
\end{eqnarray}

We must verify that the statement is true for $k$, i.e.

\begin{eqnarray}
\label{kderivada}
D^k \Theta (\psi )(\beta _1, \ldots , \beta _k)= \Theta (\psi ) ((\;.\;)\beta _1 .\ldots . \beta _k)\;,\;\; \beta _1 ,\ldots, \beta_k \in C^{\gamma}\,.
\end{eqnarray}

By the induction hypothesis \eqref{k-1-derivada}, given $\beta _1, \ldots \beta _{k-1}, \beta _k$ and $\varphi $ in $C^{\gamma}$ we have,

$D^{k-1} \Theta (\psi +\beta _k)(\beta _1, \ldots , \beta _{k-1})(\varphi)-D^{k-1} \Theta (\psi )(\beta _1, \ldots , \beta _{k-1})(\varphi )=\\
\Theta (\psi +\beta _k) (\varphi\beta _1 \ldots  \beta _{k-1})-\Theta (\psi ) (\varphi\beta _1 \ldots \beta _{k-1})$

Thus, in view of equation~(\ref{representacaoemseriesdepotencias}) we have

$$D^{k-1} \Theta (\psi +\beta _k)(\beta _1, \ldots , \beta _{k-1})(\varphi)-D^{k-1} \Theta (\psi )(\beta _1, \ldots , \beta _{k-1})(\varphi )=\sum \limits_{n=1} \limits^{\infty}\frac{1}{n!}\Theta (\psi)(\varphi \beta _1 \ldots \beta _{k-1} [\beta _k]^n)
$$

Clearly, the above equation show that

$$D^{k-1} \Theta (\psi +\beta _k)(\beta _1, \ldots , \beta _{k-1})-D^{k-1} \Theta (\psi )(\beta _1, \ldots , \beta _{k-1})= \Theta (\psi )(\;.\;)\beta _1 \ldots \beta _{k-1}\beta _k)+E_k(\beta _k)(\beta _1 \ldots \beta _{k-1})$$

Where $E_k(\beta _k)$ is the element of $L^{k-1}$ given by,

$$E_k(\beta _k)(\beta _1 \ldots \beta _{k-1})=\sum \limits_{n=2} \limits^{\infty}\frac{1}{n!}\Theta (\psi)((\;.\;) \beta _1 \ldots \beta _{k-1} [\beta _k]^n)
$$

We can use inequality~(\ref{dominando}), i.e.,
$$\mid \Theta (\psi )((\;.\;) \beta _1 \ldots \beta _k ) \mid _V \leq \mid {\mathscr L}_{\psi} \mid _V 2^k \mid \beta _1 \mid _{\gamma}\ldots \mid \beta _k \mid _{\gamma}
$$
to conclude that the map $( \beta _1 \ldots \beta _k )\mapsto  \Theta (\psi )((\;.\;) \beta _1 \ldots \beta _k )$ is an element of $L^k$.

It also follows from inequality~(\ref{dominando}) and the definition of $ E_k(\beta _k)$ that

$\mid E_k(\beta _k)(\beta _1 \ldots \beta _{k-1}) \mid _V \leq \sum \limits_{n=2}\limits^{\infty} \frac{1}{n!}\mid \mathscr L _{\psi}\mid _V 2^{k-1} \mid \beta _1 \mid _{\gamma} \ldots \mid \beta _{k-1} \mid _{\gamma}(2\mid \beta _k\mid _{\gamma})^n$\,.

Using $ \sum \limits_{n=2}\limits^{\infty} \frac{1}{n!} (2\mid \beta _k \mid _{\gamma})^n =e^{2\mid \beta _k \mid _{\gamma}}-2\mid \beta _k \mid _{\gamma}-1$ we may conclude that
$$
\mid E_k(\beta _k) \mid _{L^k}\leq \mid \mathscr L _{\psi}\mid _V 2^{k-1}\left\{ e^{2\mid \beta _k \mid _{\gamma}}-2\mid \beta _k \mid _{\gamma}-1  \right\}\,.
$$

Hence, $\frac{\mid E_k (\beta _k)\mid _V}{\mid \beta _k \mid _{\gamma}}\longrightarrow 0$, if $\mid \beta _k \mid _{\gamma} \longrightarrow 0$. Therefore, $D^k \Theta (\psi)(\beta _1, \ldots \beta _k)=\Theta (\psi )((\;.\;)\beta _1 \ldots \beta _k)$.

This ends the proof of the statement, which allow us to finish the proof of Theorem~(\ref{principal})

\qed

Let us now present a Taylor expansion of order 2 of the map $\psi \rar \Theta(\psi)$, in an integral form:

From \eqref{representacaoemseriesdepotencias} and \eqref{Taylor} we see that

$$\mathscr L _{\psi+\beta} (\varphi)(x)=
\sum \limits_{n=0} \limits^{\infty}\frac{1}{n!}\int \limits_M e^{\psi (ax)}\varphi (ax)[\beta (ax)]^n d\mu (a)\,.
$$
If we set $o_3(\beta)= \sum \limits_{n=3} \limits^{\infty}\frac{1}{n!}\int \limits_M e^{\psi (ax)}\varphi (ax)[\beta (ax)]^n d\mu (a)$ we know that $\lim_{\beta \rar 0} \frac{|o_3(\beta)|_V}{|\beta|_{\gamma}^3}=0$ and we can write

\begin{equation}\label{pse}
    \mathscr L _{\psi+\beta} (\varphi)(x)= \int e^{\psi(ax)}\varphi(ax) \left( 1+\beta(ax)+\frac{\beta(ax)^2}{2}\right)d\mu(a) + o_3(\beta)\,.
\end{equation}


\bigskip

\section{The Generalized Ruelle's Operator}\label{subshift}



In what follows, it will be convenient to adopt the following notation, $x(k)=x_k,\;x\in {\cal B},\; k \in \{0,1,2,\ldots\}.$

We shall denote by $\sigma:\cal B \rar \cal B$ the shift on $\cal B$, i.e., $\sigma (x)(k)=x(k+1),\; \forall k\geq 0 $.

Since $d_c (\sigma (x),\sigma (y))\leq d_c (x,y)$ we have that $\sigma :{\cal B}\hookleftarrow$ is a continuous operator.

\begin{defnc}
\label{seccionalmentetrivial}
Let $A:M\times M \rightarrow \mathbb R$ be a continuous function and $I$ a closed subset of $\mathbb R$.

\begin{enumerate}
\item ${\cal B}(A,I)$ is the set given by
$$
{\cal B}(A,I)=\left\{ x\in {\cal B}\; ; \;\;A(x(k),x(k+1)) \in I,\;k=0,1,2\ldots\right\}.
$$

\item For any $m \in M,\;s(m)$ the section of $m$ in $A^{-1}(I)$ is the subset of $M$ given by,

$$s(m)=\left\{a\in M\;;\;\; A(a,m)\in I \right\}\,.$$

\item  A is said sectionally trivial on $I$, when $s(m)=s(m')$ for every $m$ and $m'$ sufficiently close $m,m'\in M$, where $s(m)$ is the section of $m$ in $A^{-1}(I)$.
\end{enumerate}

\end{defnc}

Clearly, $s(m)$ is a closed subset of M.
The section of $m$ can be interpreted as  the set of pre-images of a symbol $m$ by the shift map, under restrictions on the transitions given by $A$ and $I$.

\noindent {\bf{Examples.}}


\begin{enumerate}
\label{exemplos}

\item  For $A=1$ and $I=\mathbb R$ we have ${\cal B}(A,I)={\cal B}, s(m)=M$, for all $m\in M$. Thus, obviously $A$ is sectionally trivial.

\item   Now let us consider the particular case of $M=\left\{1,\ldots , k \right\}$ equipped with the discrete distance $d$ given by
\[ d(n,j)=  \left \{   \begin{array}{cccc}
1, \mbox{ if }  n=j\,,\\
0, \mbox{ if } n\neq j\,.\end{array}
\right .    \]

Naturally, any function $A:M\times M\rightarrow \left\{ 0,1 \right\}$ is a continuous function. Of course, A can be thought of as a square matrix of order k with their entries in the set $\left\{ 0,1 \right\}$.

If we choose $I= \left\{ 1 \right\}$ then we have
$$
{\cal B}(A,I)=\left\{ x\in {\cal B}\;\;\; A(x(k),x(k+1))=1 \right\}\,.$$

We note that $s(j)=s(n)$ if $d(j,n) <1$, where $s(j)$ is the section of $j$ in $A^{-1}(I)$, since $j = n$, if $d(n,j) <1$. Therefore, A is sectionally trivial.

\item Now we will present an example which was introduced in \cite{LMMS}. The main motivation is to present a model for thermodynamic formalism for an alphabet given by $\mathbb{Z}$ or $\mathbb{N}$.
    See \cite{Sar2} for a good introduction on the subject.

Suppose $M_0 = \{z_i, i \in \mathbb{N} \}$ is an infinite sequence of points in $[0,1)$,
and suppose $\lim z_i = 1 \equiv z_{\infty}$.
Each point of $M_0$ is isolated, and there is only one accumulating point $z_{\infty}=1$. 
Then $M=M_0 \cup \{ 1\}$, with the euclidean metric, is a compact set.
Fixed a sequence $p_1,p_2,...$ such that $\sum_{i \in \mathbb{N}} p_i =1 $ and $p_i \geq 0$.
We can define the a-priori probability measure on $\ber = M^{\mathbb{N}}$ given by $\nu = \sum_{i \in \mathbb{N}} p_i. \delta_{z_i}$.
Therefore, we have a state space $M_0$ that can be identified with $\mathbb{N}$, and $M$ has a special point $z_{\infty}=1$ playing the role of the infinity.


In \cite{LMMS}, it is proved that if $\psi:\ber_0 \to\mathbb{R}$ is a Holder continuous potential,
where $\ber_0 = M_0^{\mathbb{N}}$ ,
then it can be extended to a Holder continuous function $\psi:\ber \to \mathbb{R}$. Also,
 there exists maximizing measures $\mu$ for $\psi$, and under mild assumptions it is proved in \cite{LMMS} that any maximizing measure for $\psi$ has support on $\ber_0$. Therefore the compactification of $\mathbb{N}$ is used to prove the existence of maximizing measures, and after it is proved that those measures are indeed supported on $\mathbb{N}$.

If we suppose $M_0 = \{z_i, i \in \mathbb{N} \}$ is an infinite sequence of points in $S^1\verb"\"z_{\infty}$, where $z_{\infty}$ is a single point of $S^1$ playing the role of infinity, then we have a model that can be identified with $\mathbb{Z}$ instead of $\mathbb{N}$.

We can also suppose, as a generalization of the example above, that the support of the a-priori measure depends on the point $x \in M_0$,
which means the a-priori measure in the point $x$ is given by
 $\nu_x = \sum_{i \in \mathbb{N}}
 p^x_i. \delta_{z_i}$, where $p^x_i = p_i$ if $A(p_i,x)=1$ or $p^x_i=0$ if $A(p_i,x)=0$.
If we suppose  the support of $\nu_x$ is the same for $x$ large enough, then it corresponds to define a sectionally trivial function $A:M \times M \to \mathbb{R}$.

Related results on Ergodic Optimization and thermodynamic formalism for an alphabet given by $\mathbb{Z}$ or $\mathbb{N}$ can be found on the works of Sarig (for example \cite{Sar2}). Other works are \cite{BGar,Jenkinson,JMUrb} and references therein.



\end{enumerate}

\begin{lema}
\label{iaecompacto}

 If $A:M\times M\rightarrow \mathbb R$ is a continuous function and $I$ is a closed subset of the real line then ${\cal B}(A,I)$ is a compact metric space.

\end{lema}

\noindent {\bf{Proof.}} In fact, since ${\cal B}$ is compact, it is enough to show that ${\cal B}(A,I)$ is a closed subset of ${\cal B}$. For that, consider $\left\{ x_n \right\}$ in ${\cal B}(A,I)$, with $x_n \longrightarrow x,\; x\in {\cal B}$. All we have to do is verifiy that $x$ lies in ${\cal B}(A,I)$.

From $x_n \longrightarrow x$ it follows that $x_n (k) \longrightarrow x(k),\;k=0,1,2,\ldots $.

On the other hand, for any fixed $k \in \mathbb N$ we have
$$
(x_n (k),x_n (k+1)) \in I\;,\;\; n=0,1,2,\ldots $$

When $n$ goes to infinity, we obtain
$$
A(x (k),x (k+1)) \in I $$
since $A$ is a continuous function and $I$ is a closed subset of the real line.

Using that $k$ is arbitrary, we have $x\in {\cal B}(A,I)$

\qed

Unless otherwise explicitly stated, any function $A:M\times M \rightarrow \mathbb R$ in what follows is sectionally trivial on $I$, where $I$ is a certain closed subset of the real line.

We define $\dot{C}^{\gamma}=C^{\gamma}({\cal B}(A,I),\mathbb R)$ entirely analogous to the definition of $C^{\gamma}=C({\cal B},\mathbb R)$ but with ${\cal B}(A,I)$ instead of  ${\cal B}$. By abuse of notation, the norm of $\dot{C}^{\gamma}$ will be also denoted by $\mid . \mid _{\gamma }$.

\begin{obs}
\label{holdernormdot}
Naturally, the properties of Remark~(\ref{holdernorm}) are also true in the present context:
\begin{enumerate}
\item $\mid \varphi \psi \mid _{\gamma}\leq  2 \mid \varphi \mid _{\gamma} \mid \psi \mid _{\gamma}$, for all $\varphi$ and $\psi$ in ${\dot{C}}^{\gamma}$.
\item $\mid \beta _1  \ldots \beta _k \mid _{\gamma} \leq 2^{k-1} \mid \beta _1 \mid _{\gamma} \ldots  \mid \beta _k \mid _{\gamma}$, for all $ \beta _1, \ldots , \beta _k$ in ${\dot{C}}^{\gamma}$.

\end{enumerate}

\end{obs}

\begin{defnc}
\label{sruelle}
For any $\psi \in {\dot{C}}^{\gamma}, \; 0\leq \gamma \leq 1$, the Generalized Ruelle's Operator ${\mathscr{L}}_{\psi}:{\dot{C}}^{\gamma}\hookleftarrow$ associated to $\psi$ is given by
$$
{\mathscr{L}}_{\psi}(\varphi )(x)=\int \limits_{s(x(0))} e^{\psi (ax)} \varphi (ax) d\mu (a)\;,\;\;\varphi \in {\dot{C}}^{\gamma},\; x \in {\cal B}(A,I)
$$
where $s(x(0))$ is the section of $x(0)$ in $A^{-1}(I)$ and  $ax=(a,x_0,x_1,\ldots),\;x=(x_0,x_1,\ldots )\in {\cal B}(A,I),\; a \in s(x(0))$.

\end{defnc}

\begin{lema}
\label{papeldasecao}
We summarize below some properties of the section $s(x(0)),\;x\in {\cal B}(A,I)$.

\begin{enumerate}

\item  For all $x\in {\cal B}(A,I)$, we have that

$$ax \in {\cal B}(A,I)\mbox{ if and only if } a\in s(x(0))$$

\item  If $\dot{\sigma }:{\cal B}(A,I) \hookleftarrow$ is given by
$\dot{\sigma } (x) = \sigma (x),x\in {\cal B}(A,I)$  then ${\dot{\sigma }}^{-1} ( \left\{ x\right\} )=\left\{ ax\; ; \;\; a \in s(x(0))\right\}$.

\item $s(x(0))=s(y(0))$ if $x$ and $y$ are sufficiently close, $x,y \in {\cal B}(A,I)$.

\end{enumerate}

\end{lema}

\noindent {\bf{Proof.}} we shall provide the proof of each item above.

\noindent {\it{1.}} clearly $ax \in {\cal B}(A,I)$ if and only if $A(a,x(0))\in I$, since $x\in {\cal B}(A,I)$. Thus, $ax\in {\cal B}(A,I)$ if and only if $a\in s(x(0))$.

\noindent {\it{2.}} Certainly, ${\dot{\sigma }}^{-1} ( \left\{ x\right\} )={\sigma }^{-1}( \left\{ x\right\} ) \cap {\cal B}(A,I)$ and by  ({\it{1}} ) above ${\sigma }^{-1}( \left\{ x\right\} ) \cap {\cal B}(A,I)=\left\{ ax\; ; \;\; a \in s(x(0))\right\}$. Thus,  ${\dot{\sigma }}^{-1} ( \left\{ x\right\} )=\left\{ ax\; ; \;\; a \in s(x(0))\right\}$.

\noindent {\it{3.}} From $d(x(0),y(0) \leq d_c(x,y)$ it follows that $x(0)$ and  $y(0)$ are sufficiently close, provided that $x$ and  $y$ are sufficiently close. Since $A$ is sectionally trivial the result follows

\qed

\begin{obs}
\label{sbemdefinida}
For any $\psi \in {\dot{C}}^{\gamma}, \; 0\leq \gamma \leq 1$, the Generalized Ruelle's Operator ${\mathscr{L}}_{\psi}:{\dot{C}}^{\gamma}\hookleftarrow$ associated to $\psi$ is well defined and it is a bounded linear operator. More precisely, we have,

\begin{enumerate}

\item  If $\psi \in {\dot{C}}^{0}$ then ${\mathscr{L}}_{\psi}(\varphi ) \in{\dot{C}}^{0}$, for all $\varphi \in{\dot{C}}^{0}$.

\item If $\psi \in {\dot{C}}^{\gamma},\; 0<\gamma \leq 1$, then ${\mathscr{L}}_{\psi}(\varphi ) \in {\dot{C}}^{\gamma}$, for all $\varphi \in {\dot{C}}^{\gamma}$.

\item  If $\psi \in {\dot{C}}^{}$ then ${\mathscr{L}}_{\psi}:{\dot{C}}^{0}\hookleftarrow$ is a bounded linear operator.

\item  If $\psi \in {\dot{C}}^{\gamma},\; 0<\gamma \leq 1$, then  ${\mathscr{L}}_{\psi}:{\dot{C}}^{\gamma}\hookleftarrow$ is a bounded linear operator.

\end{enumerate}

\end{obs}

\noindent {\bf{Proof.}} The arguments involved are quite similar to those of Remark~(\ref{bemdefinida}).

\noindent {\it{1.}} By ({\it{3}}) of Lemma~(\ref{papeldasecao}) for every $x$ and $y$ sufficiently close, $x,y \in {\cal B}(A,I)$ we have $s(x(0))=s(y(0))$. Thus, for every $\psi, \varphi \in {\dot{C}}^{0}$ and $x,y \in {\cal B}(A,I)$ sufficiently close, we have,

$|{\mathscr{L}}_{\psi}(\varphi )(x)-{\mathscr{L}}_{\psi}(\varphi )(y)|=|\int \limits_{s(x(0))}e^{\psi (ax)} \varphi (ax) d\mu (a)- \int \limits_{s(x(0))} e^{\psi (ay)} \varphi (ay) d\mu (a)|\leq \\
 \int \limits_{s(x(0))}|e^{\psi (ax)} \varphi (ax) - e^{\psi (ay)} \varphi (ay)| d\mu (a)$

 From now on, we continue exactly as in the proof of ({\it{1}}) of Remark~(\ref{bemdefinida}).

\noindent {\it{2.}} By ({\it{1}}) above for every $\psi, \varphi \in {\dot{C}}^{0}$ and $x,y \in {\cal B}(A,I)$ sufficiently close, we have,

$|{\mathscr{L}}_{\psi}(\varphi )(x)-{\mathscr{L}}_{\psi}(\varphi )(y)| \leq
 \int \limits_{s(x(0))}|e^{\psi (ax)} \varphi (ax) - e^{\psi (ay)} \varphi (ay)| d\mu (a)$

 Hereafter, we continue like in the proof of ({\it{2}}) of Remark~(\ref{bemdefinida}).

\noindent {\it{3.}} The proof of ({\it{3}}) of Remark~(\ref{bemdefinida}) does apply immediately.

\noindent {\it{4.}} The proof of ({\it{3}}) of Remark~(\ref{bemdefinida}) also works in our current context.

In the sequel, we shall denote by $\dot{V}$ the set of all continuous linear transformations $l:{\dot{C}}^{\gamma}\hookleftarrow$, equipped with the usual norm,

$$\mid l \mid _{\dot{V}} =\sup \limits_{\mid \varphi \mid _{\gamma}\leq 1} \mid l( \varphi )\mid _{\gamma}$$

In the present context we have the following version of Theorem~(\ref{principal}),

\begin{teo}
\label{sprincipal}

The function $\dot{\Theta } :{\dot{C}}^{\gamma}\rightarrow \dot{V}$, given by $\dot{\Theta } (\psi )={\mathscr{L}}_{\psi}$ is an analytic function. Where ${\mathscr{L}}_{\psi} $
 is the Generalized Ruelle's Operator associated to $\psi$.
\end{teo}

\noindent {\bf{Proof.}} In view of Remark~(\ref{holdernormdot}) and Remark~(\ref{sbemdefinida}) it is an easy adaptation of the proof of Theorem~(\ref{principal}).

\qed

Consider $M=(M,d)$ given by example (2), i.e.,  $M=\left\{1,\ldots , k \right\}$ and $d$ is the discrete distance. Let $\mu$ be the counting measure on $M$. For this particular case, Theorem~(\ref{sprincipal}) can be used to obtain, as a particular case, the proof of analiticity of the Ruelle operator, which was originally proved pointed out in   \cite{ma} (see also \cite{Sil}).

\begin{flushleft}

{\sc Raderson Rodrigues da Silva\\
Departamento de Matematica\\
Universidade de Brasilia\\
Campus Universitario Darcy Ribeiro-Asa Norte\\
70910-900  Brasilia DF\\
Brazil}\\
{\it email rad@mat.unb.br}

\bigskip

\bigskip

{\sc Eduardo Ant\^onio da Silva\\
Departamento de Matematica\\
Universidade de Brasilia\\
Campus Universitario Darcy Ribeiro-Asa Norte\\
70910-900  Brasilia DF\\
Brazil}\\
{\it email eduardo23maf@gmail.com}

\bigskip

\bigskip

{\sc Rafael Rigão Souza\\
Departamento de Matematica\\
Universidade Federal do Rio Grande do Sul - UFRGS\\
Avenida Bento Gonçalves, 9500\\
91509-900 - Porto Alegre RS\\
Brazil}\\
{\it email rafars@mat.ufrgs.br}

\end{flushleft}

\end{document}